\documentclass{article}
\usepackage{amssymb}
\usepackage{amsfonts}
\usepackage{amsmath}

\begin{document}

\begin{center}
{\Large Quaternion algebras and the generalized Fibonacci-Lucas quaternions}

\begin{equation*}
\end{equation*}%
Cristina FLAUT and Diana SAVIN 
\begin{equation*}
\end{equation*}
\end{center}

\textbf{Abstract. }{\small In this paper, we introduce the generalized
Fibonacci-Lucas quaternions and we prove that the set of these elements is
an order--in the sense of ring theory-- of a quaternion algebra. Moreover,
we investigate some properties of these elements.}


\medskip

\begin{equation*}
\end{equation*}

\bigskip

\textbf{1. Preliminaries}%
\begin{equation*}
\end{equation*}%
Let $\left( f_{n}\right) _{n\geq 0}$ be the Fibonacci sequence: 
\begin{equation*}
f_{0}=0;f_{1}=1;f_{n}=f_{n-1}+f_{n-2},\;n\geq 2
\end{equation*}%
and $\left( l_{n}\right) _{n\geq 0}$ be the Lucas sequence: 
\begin{equation*}
l_{0}=2;f_{1}=1;l_{n}=l_{n-1}+l_{n-2},\;n\geq 2
\end{equation*}%
In the paper [Ho; 61], A. F. Horadam generalized Fibonacci numbers in the
following way: 
\begin{equation*}
h_{0}=p,h_{1}=q,h_{n}=h_{n-1}+h_{n-2},\;n\geq 2,
\end{equation*}%
where $p$ and $q$ are arbitrary integer numbers.\newline
In the paper [Ho; 63], A. F. Horadam introduced the Fibonacci quaternions
and generalized Fibonacci quaternions. In their work [Fl, Sh; 12], C. Flaut
and V. Shpakivskyi and later in the paper [Ak, Ko, To; 14], M. Akyigit, H. H
Kosal, M. Tosun found some properties of the generalized Fibonacci
quaternions. In the papers [Fl, Sa; 14] and [Fl, Sa, Io; 13], the authors
gave the definitions of the Fibonacci symbol elements and Lucas symbol
elements and they determined many properties of these elements. All these
elements determine a $\mathbb{Z}$-module structure.\newline
In this paper, we define the generalized Fibonacci-Lucas numbers and the
generalized Fibonacci-Lucas quaternions and we prove that, in the second
case, they determine an order of a quaternion algebra.

For other results regardind above notions, the reader is referred to [Fl;
06], [Sa, Fl, Ci; 09], [Fl, Sh; 13(1)], [Fl, Sa; 14], [Fl, Sa; 15],[Fl, St;
09]. 
\begin{equation*}
\end{equation*}

\textbf{2. Properties of the Fibonacci and Lucas numbers}

\begin{equation*}
\end{equation*}

From [Fib.], the following properties of Fibonacci numbers are known:

\medskip

\textbf{Proposition 2.1.} \textit{Let} $(f_{n})_{n\geq 0}$ \textit{be the
Fibonacci sequence} \textit{and let } $(l_{n})_{n\geq 0}$ \textit{be the
Lucas sequence.} \textit{Therefore the following properties hold:}\newline
i) 
\begin{equation*}
f_{n}^{2}+f_{n+1}^{2}=f_{2n+1},\forall ~n\in \mathbb{N};
\end{equation*}%
ii) 
\begin{equation*}
f_{n+1}^{2}-f_{n-1}^{2}=f_{2n},\forall ~n\in \mathbb{N}^{\ast };
\end{equation*}%
iii) 
\begin{equation*}
l_{n}^{2}-f_{n}^{2}=4f_{n-1}f_{n+1},\forall ~n\in \mathbb{N}^{\ast };
\end{equation*}%
iv) 
\begin{equation*}
l_{n}^{2}+l_{n+1}^{2}=5f_{2n+1},\forall ~n\in \mathbb{N};
\end{equation*}%
v) 
\begin{equation*}
l_{n}^{2}=l_{2n}+2\left( -1\right) ^{n},\forall ~n\in \mathbb{N}^{\ast };
\end{equation*}%
vi) 
\begin{equation*}
f_{n+1}+f_{n-1}=l_{n},\forall ~n\in \mathbb{N}^{\ast };
\end{equation*}%
vii) 
\begin{equation*}
l_{n}+l_{n+2}=5f_{n+1},\forall ~n\in \mathbb{N};
\end{equation*}%
viii) 
\begin{equation*}
f_{n}+f_{n+4}=3f_{n+2},\forall ~n\in \mathbb{N};
\end{equation*}%
ix) 
\begin{equation*}
f_{m}l_{m+p}=f_{2m+p}+\left( -1\right) ^{m+1}\cdot f_{p},\forall ~m,p\in 
\mathbb{N};
\end{equation*}%
x) 
\begin{equation*}
f_{m+p}l_{m}=f_{2m+p}+\left( -1\right) ^{m}\cdot f_{p},\forall ~m,p\in 
\mathbb{N};
\end{equation*}%
xi) 
\begin{equation*}
f_{m}f_{m+p}=\frac{1}{5}\left( l_{2m+p}+\left( -1\right) ^{m+1}\cdot
l_{p}\right) ,\forall ~m,p\in \mathbb{N};
\end{equation*}%
xii) 
\begin{equation*}
l_{m}l_{p}+5f_{m}f_{p}=2l_{m+p},\forall ~m,p\in \mathbb{N}.
\end{equation*}%
xiii) 
\begin{equation*}
f_{n}\ are\ even\ numbers\ if\ and\ only\ if\ n\equiv 0\ (mod\ 3).
\end{equation*}

\medskip

In the next proposition, we will give other elementary properties of the
Fibonacci and Lucas numbers. These properties will be necessary in the
following.\medskip \newline
\textbf{Proposition 2.2.} \textit{Let} $(f_{n})_{n\geq 0}$ \textit{be the
Fibonacci sequence} \textit{and} $(l_{n})_{n\geq 0}$ \textit{be the Lucas
sequence} \textit{Then we have that}\newline
\begin{equation*}
l_{m}l_{m+p}=l_{2m+p}+\left( -1\right) ^{m}\cdot l_{p},\forall ~m,p\in 
\mathbb{N}.
\end{equation*}%
\textbf{Proof.} If we denote $\alpha =\frac{1+\sqrt{5}}{2}$ and $\beta =%
\frac{1-\sqrt{5}}{2},$ using Binet's formula, we have $f_{n}=\frac{\alpha
^{n}-\beta ^{n}}{\alpha -\beta }=\frac{1}{\sqrt{5}}\left( \alpha ^{n}-\beta
^{n}\right) ,$ $\forall ~n\in \mathbb{N}$ and $l_{n}=\alpha ^{n}+\beta ^{n},$
$\forall ~n\in \mathbb{N}.$\newline
i) Let $m,p\in \mathbb{N},$ $p\leq m,$ therefore we have 
\begin{equation*}
l_{m}l_{p}-5f_{m}f_{p}=\left( \alpha ^{m}+\beta ^{m}\right) \cdot \left(
\alpha ^{p}+\beta ^{p}\right) -\left( \alpha ^{m}-\beta ^{m}\right) \cdot
\left( \alpha ^{p}-\beta ^{p}\right) =
\end{equation*}%
\begin{equation*}
=2\alpha ^{p}\beta ^{p}\left( \alpha ^{m-p}+\beta ^{m-p}\right) =2\left(
-1\right) ^{p}\cdot l_{m-p}.
\end{equation*}%
It results $l_{m}l_{p}-5f_{m}f_{p}=2\left( -1\right) ^{p}\cdot l_{m-p}.$
Adding this equality with the equality from Proposition 2.1 (xii)) , we
obtain $l_{m}l_{p}=l_{m+p}+\left( -1\right) ^{p}\cdot l_{m-p},\forall
~m,p\in \mathbb{N},p\leq m.$ From here, it results that $%
l_{m}l_{m+p}=l_{2m+p}+\left( -1\right) ^{m}\cdot l_{p},\forall ~m,p\in 
\mathbb{N}.$$\square \medskip $ 
\begin{equation*}
\end{equation*}%
\textbf{3. Generalized Fibonacci-Lucas numbers and generalized
Fibonacci-Lucas quaternions}%
\begin{equation*}
\end{equation*}

Let $n$ be an arbitrary positive integer and $p,q$ be two arbitrary
integers. The sequence $g_{n}$ $(n\geq 1),$ where 
\begin{equation*}
g_{n+1}=pf_{n}+ql_{n+1},\;n\geq 0
\end{equation*}%
is called \textit{the generalized Fibonacci-Lucas numbers}.\newline
To emphasize the integer $p$ and $q,$ in the following we will use the
notation $g_{n}^{p,q}$ instead of $g_{n}.$ \smallskip

Let $\mathbb{H}\left( \alpha ,\beta \right) $ be the generalized real\
quaternion algebra, the algebra of the elements of the form $a=a_{1}\cdot
1+a_{2}i+a_{3}j+a_{4}k,$ where $a_{i}\in \mathbb{R},i^{2}= \alpha ,j^{2}=
\beta ,$ $k=ij=-ji.$ We denote by $\mathbf{t}\left( a\right) $ and $\mathbf{n%
}\left( a\right) $ the trace and the norm of a real quaternion $a.$ The norm
of a generalized quaternion has the following expression $\mathbf{n}\left(
a\right) =a_{1}^{2}-\alpha a_{2}^{2}-\beta a_{3}^{2}+\alpha \beta a_{4}^{2}$
and $\mathbf{t}\left( a\right) =2a_{1}.$ It is known that for $a\in $ $%
\mathbb{H}\left( \alpha ,\beta \right) ,$ we have $a^{2}-\mathbf{t}\left(
a\right) a+\mathbf{n}\left( a\right) =0.$ The quaternion algebra $\mathbb{H}%
\left( \alpha ,\beta \right) $ is a \textit{division algebra} if for all $%
a\in \mathbb{H}\left( \alpha ,\beta \right) ,$ $a\neq 0$ we have $\mathbf{n}%
\left( a\right) \neq 0,$ otherwise $\mathbb{H}\left( \alpha ,\beta \right) $
is called a \textit{split algebra}.

Let $\mathbb{H}_{\mathbb{Q}}\left( \alpha ,\beta \right) $ be the
generalized quaternion algebra over the rational field$.$ We define the $n$-%
\textit{th generalized Fibonacci-Lucas quaternion} to be the element of the
form 
\begin{equation*}
G_{n}^{p,q}=g_{n}^{p,q}\cdot 1+g_{n+1}^{p,q}\cdot i+g_{n+2}^{p,q}\cdot
j+g_{n+3}^{p,q}\cdot k,
\end{equation*}%
where $i^{2}=\alpha ,j^{2}=\beta ,$ $k=ij=-ji.$ \smallskip

In the following proposition, for $\alpha =-1$ and $\beta =p,$ we compute
the norm for the $n$-th generalized Fibonacci-Lucas quaternions.\medskip

\textbf{Proposition 3.1.} \textit{Let} $n,p$ \textit{be two arbitrary
positive integers and} $q$ \textit{be an arbitrary integer. Let }$%
G_{n}^{p,q} $ \textit{be the} $n$-\textit{th generalized Fibonacci-Lucas
quaternion. Then the norm of the element} $G_{n}^{p,q}$ \textit{in the
quaternion algebra} $\mathbb{H}_{\mathbb{Q}}\left( -1,p\right) $ \textit{has
the form}\newline
i) 
\begin{equation*}
\mathbf{n}\left( G_{n}^{p,q}\right) =\left(
p^{3}+p^{2}+8p^{2}q+15pq^{2}+2pq\right) f_{2n-1}+
\end{equation*}%
\begin{equation*}
+\left( -3p^{3}+5q^{2}-22p^{2}q-40pq^{2}+2pq\right) f_{2n+1};
\end{equation*}%
\textit{or\newline
}ii) 
\begin{equation*}
\mathbf{n}\left( G_{n}^{p,q}\right) =g_{2n}^{a,b},
\end{equation*}%
\textit{where} $a=4p^{3}+p^{2}+30p^{2}q+55pq^{2}-5q^{2}$ \textit{and} $%
b=-3p^{3}+5q^{2}-22p^{2}q-40pq^{2}+2pq.$$\medskip $

\textbf{Proof.} i) 
\begin{equation*}
\mathbf{n}\left( G_{n}^{p,q}\right) =G_{n}^{p,q}\overline{G_{n}^{p,q}}%
=\left( g_{n}^{p,q}\right) ^{2}+\left( g_{n+1}^{p,q}\right) ^{2}-p\left(
g_{n+2}^{p,q}\right) ^{2}-p\left( g_{n+3}^{p,q}\right) ^{2}=
\end{equation*}%
\begin{equation*}
=\left( pf_{n-1}+ql_{n}\right) ^{2}+\left( pf_{n}+ql_{n+1}\right)
^{2}-p\left( pf_{n+1}+ql_{n+2}\right) ^{2}-p\left( pf_{n+2}+ql_{n+3}\right)
^{2}=
\end{equation*}%
\begin{equation*}
=p^{2}\left( f_{n-1}^{2}+f_{n}^{2}\right) -p^{3}\left(
f_{n+1}^{2}+f_{n+2}^{2}\right) +q^{2}\left( l_{n}^{2}+l_{n+1}^{2}\right)
-pq^{2}\left( l_{n+2}^{2}+l_{n+3}^{2}\right) +
\end{equation*}%
\begin{equation*}
+2pq\left( f_{n-1}l_{n}+f_{n}l_{n+1}\right) -2p^{2}q\left(
f_{n+1}l_{n+2}+f_{n+2}l_{n+3}\right) \eqno(3.1).
\end{equation*}%
Using Proposition 2.1 (i; vi;viii; ix) and the relation $(3.1),$ we obtain 
\begin{equation*}
\mathbf{n}\left( G_{n}^{p,q}\right) =p^{2}f_{2n-1}+5q^{2}f_{2n+1}+2pq\left[
f_{2n-1}+\left( -1\right) ^{n}+f_{2n+1}+\left( -1\right) ^{n+1}\right] -
\end{equation*}%
\begin{equation*}
-p^{3}f_{2n+3}-5pq^{2}f_{2n+5}-2p^{2}q\left[ f_{2n+3}+\left( -1\right)
^{n+2}+f_{2n+5}+\left( -1\right) ^{n+3}\right] =
\end{equation*}%
\begin{equation*}
=p^{2}f_{2n-1}+5q^{2}f_{2n+1}-p^{3}f_{2n+3}-5pq^{2}f_{2n+5}+2pql_{2n}-2p^{2}ql_{2n+4}.
\end{equation*}%
From Fibonacci recurrence, we obtain \newline
\begin{equation*}
f_{2n+3}=3f_{2n+1}-f_{2n-1},f_{2n+5}=8f_{2n+1}-3f_{2n-1}.\newline
\end{equation*}%
Using these equalities, Proposition 2.1 (vi) and Fibonacci recurrence, we
obtain $l_{2n}=f_{2n+1}+f_{2n-1}$ and $%
l_{2n+4}=f_{2n+3}+f_{2n+5}=11f_{2n+1}-4f_{2n-1}.$\newline
It results%
\begin{equation*}
\mathbf{n}\left( G_{n}^{p,q}\right) =\left(
p^{3}+p^{2}+8p^{2}q+15pq^{2}+2pq\right) f_{2n-1}+
\end{equation*}%
\begin{equation*}
+\left( -3p^{3}+5q^{2}-22p^{2}q-40pq^{2}+2pq\right) f_{2n+1}.
\end{equation*}%
ii) Using Proposition 2.1 (vi), last equality becomes 
\begin{equation*}
\mathbf{n}\left( G_{n}^{p,q}\right) =\left(
4p^{3}+p^{2}+30p^{2}q+55pq^{2}-5q^{2}\right) f_{2n-1}+
\end{equation*}%
\begin{equation*}
+\left( -3p^{3}+5q^{2}-22p^{2}q-40pq^{2}+2pq\right) l_{2n}.
\end{equation*}%
If we denote $a=4p^{3}+p^{2}+30p^{2}q+55pq^{2}-5q^{2}$ and $%
b=-3p^{3}+5q^{2}-22p^{2}q-40pq^{2}+2pq,$ the last equality is equivalent
with $\mathbf{n}\left( G_{n}^{p,q}\right)
=af_{2n-1}+bl_{2n}=g_{2n}^{a,b}.\square \medskip $

\smallskip Let $A$ be a Noetherian integral domain with the field of the
fractions $K$ and let $\mathbb{H}_{K}\left( \alpha ,\beta \right) $ be the
generalized quaternion algebra. We recall that $\mathcal{O}$ is an order in $%
\mathbb{H}_{K}\left( \alpha ,\beta \right) $ if $\mathcal{O}$ $\subseteq $ $%
\mathbb{H}_{K}\left( \alpha ,\beta \right) $ is an $A$-lattice of $\mathbb{H}%
_{K}\left( \alpha ,\beta \right) $ (i.e a finitely generated $A$ submodule
of $\mathbb{H}_{K}\left( \alpha ,\beta \right) $) which is also a subring of 
$\mathbb{H}_{K}\left( \alpha ,\beta \right) $ (see [Vo; 10]).\newline
In the following, we will built an order of a quaternion algebras using the
generalized Fibonacci-Lucas quaternions. Also we will prove that
Fibonacci-Lucas quaternions can have an algebra structure over $\mathbb{Q}$.
For this, we make the following remarks.\medskip

\textbf{Remark 3.1.} \textit{Let} $n$ \textit{be an arbitrary positive
integer and} $p,q$ \textit{be two arbitrary integers. Let }$\left(
g_{n}^{p,q}\right) _{n\geq 1}$ \textit{be the generalized Fibonacci-Lucas
numbers. Then} 
\begin{equation*}
pf_{n+1}+ql_{n}=g_{n}^{p,q}+g_{n+1}^{p,0},\forall ~n\in \mathbb{N}^{\ast }.
\end{equation*}%
\textbf{Proof.} 
\begin{equation*}
pf_{n+1}+ql_{n}=pf_{n-1}+ql_{n}+pf_{n}=g_{n}^{p,q}+pf_{n}=g_{n}^{p,q}+g_{n+1}^{p,o}.
\end{equation*}%
$\square \medskip $ \smallskip

\textbf{Remark 3.2.} \textit{Let} $n$ \textit{be an arbitrary positive
integer and} $p,q$ \textit{be two arbitrary integers. Let }$\left(
g_{n}^{p,q}\right) _{n\geq 1}$ \textit{be the generalized Fibonacci-Lucas
numbers and} $\left( G_{n}^{p,q}\right) _{n\geq 1}$ \textit{be the
generalized Fibonacci-Lucas quaternion elements. Then:} 
\begin{equation*}
G_{n}^{p,q}=0\ if\ and\ only\ if\ p=q=0.
\end{equation*}%
\textbf{Proof.} " $\Leftarrow $" It is trivial.\newline
" $\Rightarrow $" If $G_{n}^{p,q}=0,$ since $\left\{ 1,i,j,k\right\} $ is a
basis in $\mathbb{H}_{\mathbb{Q}}\left( \alpha ,\beta \right) ,$ we obtain
that $g_{n}^{p,q}=g_{n+1}^{p,q}=g_{n+2}^{p,q}=g_{n+3}^{p,q}=0.$ It results $%
g_{n-1}^{p,q}=g_{n+1}^{p,q}-g_{n}^{p,q}=0,$ ..., $g_{2}^{p,q}=0$, $%
g_{1}^{p,q}=0.$ But $g_{1}^{p,q}=pf_{0}+ql_{1}=2q,$ therefore $q=0.$ From $%
g_{2}^{p,q}=0,$ we obtain $p=0.$$\square \medskip $

\textbf{Theorem 3.1.} \textit{Let} $M$ \ \textit{be} \textit{the set} 
\begin{equation*}
M=\left\{ \sum\limits_{i=1}^{n}5G_{n_{i}}^{p_{i},q_{i}}|n\in \mathbb{N}%
^{\ast },p_{i},q_{i}\in \mathbb{Z},(\forall )i=\overline{1,n}\right\} \cup
\left\{ 1\right\} .
\end{equation*}%
\newline
\textit{1)} \textit{The set} $M$ \textit{has a ring structure with
quaternions addition and multiplication.\newline
2)} \textit{The set} $M$ \textit{is an order of the quaternion algebra} $%
\mathbb{H}_{\mathbb{Q}}\left( \alpha ,\beta \right) .$\newline
\textit{3)} \textit{The set} $M^{\prime }=\left\{
\sum\limits_{i=1}^{n}5G_{n_{i}}^{p_{i}^{\prime },q_{i}^{\prime }} |n\in 
\mathbb{N}^{\ast },p_{i}^{\prime },q_{i}^{\prime }\in \mathbb{Q},(\forall )i=%
\overline{1,n}\right\} \cup \left\{ 1\right\} $ \textit{is a} $\mathbb{Q-}$%
\textit{algebra}.\smallskip \newline

\textbf{Proof.} 2) First, we remark that $0$$\in $$M$ (using Remark 3.2).%
\newline
Now we prove that $M$ \textit{is a} $\mathbb{Z}-$ submodule of $\mathbb{H}_{%
\mathbb{Q}}\left( \alpha ,\beta \right) .$\newline
Let $n,m\in \mathbb{N}^{\ast },$ $a,b,p,q,p^{^{\prime }},q^{^{\prime }}\in 
\mathbb{Z}.$ It is easy to prove that 
\begin{equation*}
ag_{n}^{p,q}+bg_{m}^{p^{^{\prime }},q^{^{\prime
}}}=g_{n}^{ap,aq}+g_{m}^{bp^{^{\prime }},bq^{^{\prime }}}.
\end{equation*}%
This implies that 
\begin{equation*}
aG_{n}^{p,q}+bG_{m}^{p^{^{\prime }},q^{^{\prime
}}}=G_{n}^{ap,aq}+G_{m}^{bp^{^{\prime }},bq^{^{\prime }}}.
\end{equation*}%
From here, we get immediately that $M$ \textit{is a} $\mathbb{Z}-$ submodule
of the quaternion algebra $\mathbb{H}_{\mathbb{Q}}\left( \alpha ,\beta
\right) .$ Since $\left\{ 1,i,j,k\right\} $ is a basis for this submodule,
it results that $M$ is a free $\mathbb{Z}-$ module of rank $4.$\newline
Now, we prove that $M$ is a subring of $\mathbb{H}_{\mathbb{Q}}\left( \alpha
,\beta \right) .$ It is enough to show that $5G_{n}^{p,q}\cdot
5G_{m}^{p^{^{\prime }},q^{^{\prime }}}$$\in $$M.$ For this, if $m<n,$ we
calculate 
\begin{equation*}
5g_{n}^{p,q}\cdot 5g_{m}^{p^{^{\prime }},q^{^{\prime }}}=5\left(
pf_{n-1}+ql_{n}\right) \cdot 5\left( p^{^{\prime }}f_{m-1}+q^{^{\prime
}}l_{m}\right) =
\end{equation*}%
\begin{equation*}
=25pp^{^{\prime }}f_{n-1}f_{m-1}+25pq^{^{\prime }}f_{n-1}l_{m}+25p^{^{\prime
}}qf_{m-1}l_{n}+25qq^{^{\prime }}l_{n}l_{m}\eqno(3.2).
\end{equation*}%
Using Proposition 2.1 (ix, x, xi), Proposition 2.2 , Remark 3.1 and the
equality (3.2) we obtain: 
\begin{equation*}
5g_{n}^{p,q}\cdot 5g_{m}^{p^{^{\prime }},q^{^{\prime }}}=5pp^{^{\prime }} 
\left[ l_{m+n-2}+\left( -1\right) ^{m}\cdot l_{n-m}\right] +25pq^{^{\prime }}%
\left[ f_{m+n-1}+\left( -1\right) ^{m}\cdot f_{n-m-1}\right] +
\end{equation*}%
\begin{equation*}
+25p^{^{\prime }}q\left[ f_{m+n-1}+\left( -1\right) ^{m}\cdot f_{n-m+1}%
\right] +25qq^{^{\prime }}\left[ l_{m+n}+\left( -1\right) ^{m}\cdot l_{n-m}%
\right] =
\end{equation*}%
\begin{equation*}
=5\left( pp^{^{\prime }}l_{m+n-2}+5p^{^{\prime }}qf_{m+n-1}\right) +5\left[
5p^{^{\prime }}q\left( -1\right) ^{m}\cdot f_{n-m+1}+pp^{^{\prime }}\left(
-1\right) ^{m}\cdot l_{n-m}\right] +
\end{equation*}%
\begin{equation*}
+25\left( pq^{^{\prime }}f_{m+n-1}+qq^{^{\prime }}l_{m+n}\right) +25\left[
pq^{^{\prime }}\cdot \left( -1\right) ^{m}\cdot f_{n-m-1}+qq^{^{\prime
}}\cdot \left( -1\right) ^{m}\cdot l_{n-m}\right] =
\end{equation*}%
\begin{equation*}
=5g_{m+n-2}^{5p^{^{\prime }}q,pp^{^{\prime }}}+5g_{m+n-1}^{5p^{^{\prime
}}q,0}+5g_{n-m}^{5p^{^{\prime }}q\cdot \left( -1\right) ^{m},pp^{^{\prime
}}\cdot \left( -1\right) ^{m}}+5g_{n-m+1}^{5p^{^{\prime }}q\cdot \left(
-1\right) ^{m},0}+
\end{equation*}%
\begin{equation*}
+5g_{m+n}^{5pq^{^{\prime }},5qq^{^{\prime }}}+5g_{n-m}^{5pq^{^{\prime
}}\cdot \left( -1\right) ^{m},5qq^{^{\prime }}\cdot \left( -1\right) ^{m}}.
\end{equation*}%
Therefore, we obtain that $5G_{n}^{p,q}\cdot 5G_{m}^{p^{^{\prime
}},q^{^{\prime }}}$$\in $$M.$\newline
From here, it results that $M$ is an order of the quaternion algebra $%
\mathbb{H}_{\mathbb{Q}}\left( \alpha ,\beta \right) .$

1) and 3) are obviously.$\square \medskip $

\textbf{Remark 3.3.} For $\alpha =\beta =-1,$ we have that $M$ is included
in the set of Hurwitz quaternions,%
\begin{equation*}
\mathcal{H}=\{q=a_{1}+a_{2}i+a_{3}j+a_{4}k\in \mathbb{H}_{\mathbb{Q}}\left(
-1,-1\right) ,a_{1},a_{2},a_{3},a_{4}\in \mathbb{Z}\text{ or }\mathbb{Z}+%
\frac{1}{2}\},
\end{equation*}%
which is a maximal order in $\mathbb{H}_{\mathbb{Q}}\left( -1,-1\right) .$%
\smallskip \newline

From [Lam; 04], [Sa; 14], it is know that if $p$ is an odd prime positive
integer, the algebra $\mathbb{H}_{\mathbb{Q}}\left( -1,p\right) $ is a split
algebra if and only if $p\equiv 1$ (mod $4$). In the following, we will show
that there are an infinite numbers of generalized Fibonacci-Lucas quaternion
elements which are invertible in this algebra.\smallskip \newline

\textbf{Proposition 3.2.} \textit{Let} $n$ \textit{be an arbitrary positive
integer. Let} $(f_{n})_{n\geq 0}$ \textit{be the Fibonacci sequence} \textit{%
and} $(l_{n})_{n\geq 0}$ \textit{be the Lucas sequence. Let} $p$ \textit{be
an odd prime positive integer}, $p\equiv 1$ (\textit{mod} $4$), $q$ \textit{%
be an arbitrary integer. Let }$G_{n}^{p,q}$ \textit{be the} $n$-\textit{th
generalized Fibonacci-Lucas quaternion and } $\mathbb{H}_{\mathbb{Q}}\left(
-1,p\right) $ \textit{be the quaternion algebra}$.$ \textit{The following
statements are true}:\newline
i) $\mathbf{n}\left( G_{n}^{p,q}\right) \neq 0,$ $\left( \forall \right) $ $%
\left( n,q\right) $$\in $$\mathbb{N}^{\ast }\times \mathbb{N};$\newline
ii) \textit{If} $p=5,$ \textit{then} $\mathbf{n}\left( G_{n}^{5,q}\right)
\neq 0,$ $\left( \forall \right) $ $\left( n,q\right) $$\in $$\mathbb{N}%
^{\ast }\times \mathbb{Z}_{-}$, $\left( n,q\right) \neq \left( 1,-2\right) ;$%
\newline
iii) \textit{if} $p>5$ \textit{and} $n\equiv 0$ (\textit{mod} $3$), \textit{%
then} $\mathbf{n}\left( G_{n}^{p,q}\right) \neq 0,$ $\left( \forall \right) $
$q$$\in $$\mathbb{Z}.$ $\medskip $

\textbf{Proof.} We know that an element $x\neq 0$ from a quaternion algebra
is invertible in this algebra if the norm $\mathbf{n}\left( x\right) \neq 0.$%
\newline
i) From Proposition 3.1, we know that 
\begin{equation*}
\mathbf{n}\left( G_{n}^{p,q}\right) =\left(
p^{3}+p^{2}+8p^{2}q+15pq^{2}+2pq\right) f_{2n-1}+
\end{equation*}%
\begin{equation*}
+\left( -3p^{3}+5q^{2}-22p^{2}q-40pq^{2}+2pq\right) f_{2n+1}\eqno(3.3).
\end{equation*}%
If $q\in $$\mathbb{N},$ since $p\in $$\mathbb{N}^{\ast },$ it results 
\begin{equation*}
p^{3}+p^{2}+8p^{2}q+15pq^{2}+2pq<3p^{3}-5q^{2}+22p^{2}q+40pq^{2}-2pq.
\end{equation*}%
Using the inequality $f_{2n-1}<f_{2n+1},$ we obtain that $\mathbf{n}\left(
G_{n}^{p,q}\right) <0.$ So $\mathbf{n}\left( G_{n}^{p,q}\right) \neq 0.$%
\newline
\smallskip \newline
In the case when $q\in $$\mathbb{Z}_{-},$ if there exist generalized
Fibonacci-Lucas quaternion elements with $\mathbf{n}\left(
G_{n}^{p,q}\right) =0,$ using relation $(3.3),$ we obtain that 
\begin{equation*}
p\left[ \left( p^{2}+p+8pq+15q^{2}+2q\right) f_{2n-1}+\left(
-3p^{2}-22pq-40q^{2}+2q\right) f_{2n+1}\right] =
\end{equation*}%
\begin{equation*}
=-5q^{2}f_{2n+1},\eqno(3.4).
\end{equation*}%
therefore, $p|5q^{2}f_{2n+1}.$ But we know from the hypothesis that $p$ is a
prime positive integer, $p\equiv 1$ (mod $4$), so we consider two cases: 
\newline
ii) The case~ $p=5.$ Using relation $(3.3)$ and making some calculations we
get that 
\begin{equation*}
\mathbf{n}\left( G_{n}^{5,q}\right) =0\Leftrightarrow \left(
5q^{2}+10q+30\right) f_{2n-1}-\left( 39q^{2}+108q+75\right) f_{2n+1}=0\eqno%
(3.5).
\end{equation*}%
Therefore 
\begin{equation*}
0<5q^{2}+10q+30<39q^{2}+108q+75,\ \left( \forall \right) q\in \mathbb{Z}%
,q\neq -2;-1.
\end{equation*}%
Using the fact that $f_{2n-1}<f_{2n+1},\ \left( \forall \right) n\in \mathbb{%
N}^{\ast },$ it results 
\begin{equation*}
\mathbf{n}\left( G_{n}^{5,q}\right) <0\ \left( \forall \right) n\in \mathbb{N%
}^{\ast },\ \left( \forall \right) q\in \mathbb{Z}_{-},q\neq -2;-1.
\end{equation*}%
In the following, we will find if the equations $\mathbf{n}\left(
G_{n}^{5,-1}\right) =0$ respectively $\mathbf{n}\left( G_{n}^{5,-2}\right)
=0 $ have solutions for $n\in \mathbb{N}^{\ast }.$\newline
Using the relation $(3.5)$ and the Fibonacci recurrence, we have 
\begin{equation*}
\mathbf{n}\left( G_{n}^{5,-1}\right) =0\Leftrightarrow
25f_{2n-1}-6f_{2n+1}=0\Leftrightarrow 7f_{2n-1}+6f_{2n-3}=0.
\end{equation*}%
Since $f_{2n-1},f_{2n-3}>0$ it results that does not exist $n$$\in $$\mathbb{%
N}^{\ast }$ such that $7f_{2n-1}+6f_{2n-3}=0.$\newline
Using relation $(3.5)$ and the Fibonacci recurrence, it results 
\begin{equation*}
\mathbf{n}\left( G_{n}^{5,-2}\right) =0\Leftrightarrow
30f_{2n-1}-15f_{2n+1}=0\Leftrightarrow f_{2n-2}=0\Leftrightarrow n=1.
\end{equation*}%
iii) when $p>5,$ since $n\equiv 0$ ( mod $3$) it results that $2n-1\equiv 2$
( mod $3$) and $2n+1\equiv 1$ ( mod $3$). Applying Applying Proposition 2.1
(xiii), we obtain that $f_{2n-1}$ and $f_{2n+1}$ are odd numbers. Since $p$
is an odd number, we obtain that the equality (3.4) cannot be held (is easy
to prove that the equality (3.4) holds only when $f_{2n-1}$ is an even
number or $f_{2n+1}$ is an even number).$\square \medskip $

For $a\in \mathbb{H}_{\mathbb{Q}}\left( \alpha ,\beta \right) .$ The
centralizer of the element $a$ is 
\begin{equation*}
C\left( a\right) =\{x\in \mathbb{H}_{\mathbb{Q}}\left( \alpha ,\beta \right)
~/~ax=xa\}.
\end{equation*}
From [Fl; 01], Proposition 2.13, we know that the equation 
\begin{equation}
ax=bx,a,b\in \mathbb{H}_{K}\left( \alpha ,\beta \right) ,  \tag{3.6.}
\end{equation}%
where $K$ is an arbitrary field of $charK\neq 0,$ $a,b\notin K,a\neq 
\overline{b},$ has the solutions of the form%
\begin{equation}
x\text{=}\lambda _{1}\left( a\text{-}\mathbf{t}\left( a\right) \text{+}b%
\text{-}\mathbf{t}\left( b\right) \right) \text{+}\lambda _{2}\left( \mathbf{%
n}\left( a\text{-}\mathbf{t}\left( a\right) \right) \text{-}\left( a\text{-}%
\mathbf{t}\left( a\right) \right) \left( b\text{-}\mathbf{t}\left( b\right)
\right) \right) ,\lambda _{1},\lambda _{2}\in K  \tag{3.7.}
\end{equation}%
if $\mathbb{H}_{K}\left( \alpha ,\beta \right) $ is a division quaternion
algebra or if $\mathbb{H}_{K}\left( \alpha ,\beta \right) $ is a split
quaternion algebra and $\mathbf{n}\left( a\right) \neq 0,\mathbf{n}\left(
b\right) \neq 0$.\medskip

\textbf{Proposition 3.3.} \textit{Let} $n$ \textit{be a positive integer,}$%
n\equiv 0$ (\textit{mod} $3$)\textit{. Let} $(f_{n})_{n\geq 0}$ \textit{be
the Fibonacci sequence} \textit{and} $(l_{n})_{n\geq 0}$ \textit{be the
Lucas sequence. Let} $p$ \textit{be an odd prime positive integer}, $p\equiv
1$ (\textit{mod} $4$), $q$ \textit{be an arbitrary integer. Therefore, the
centralizer of the element} $G_{n}^{p,q}\in $ $\mathbb{H}_{\mathbb{Q}}\left(
-1,p\right) $ \textit{with} $p>5,$ $q$ \textit{arbitrary integer is the set}%
\begin{equation*}
C\left( G_{n}^{p,q}\right) =\{G_{n}^{\gamma ,\delta }+\Lambda ,\Lambda \in 
\mathbb{Q}\},
\end{equation*}%
\textit{where} $\gamma =2\lambda _{1}p,\delta =2\lambda _{1}q,\Lambda
=g_{n}^{-2\lambda _{1}p,-2\lambda _{1}q}$+$g_{2n}^{2\lambda _{2}a-4\lambda
_{2}pq,2\lambda _{2}b-2\lambda _{2}q^{2}}$-$2\lambda _{2}\frac{p^{2}}{5}%
l_{2n-2}$-$2\lambda _{2}A,$ \textit{with} $\lambda _{1},\lambda _{2}\in 
\mathbb{Q}$ \textit{and} $A=\frac{2p^{2}}{5}(-1)^{n}+2q^{2}\left( -1\right)
^{n}+2pq\left( -1\right) ^{n}.\medskip $

\textbf{Proof.} Since $C\left( G_{n}^{p,q}\right) =\{x\in \mathbb{H}_{%
\mathbb{Q}}\left( -1,p\right) ~/~G_{n}^{p,q}x=xG_{n}^{p,q}\}$, using
relations $\left( 3.6\right) $ and $\left( 3.7\right) $ for $a=b,$ we obtain
that the equation $G_{n}^{p,q}x=xG_{n}^{p,q}$ has the solutions of the form 
\begin{equation}
x=2\lambda _{1}\left( G_{n}^{p,q}\text{-}\mathbf{t}\left( G_{n}^{p,q}\right)
\right) +2\lambda _{2}\mathbf{n}\left( G_{n}^{p,q}\text{-}\mathbf{t}\left(
G_{n}^{p,q}\right) \right) ,\lambda _{1},\lambda _{2}\in \mathbb{Q}. 
\tag{3.8.}
\end{equation}

For $G_{n}^{p,q}=g_{n}^{p,q}\cdot 1+g_{n+1}^{p,q}\cdot i+g_{n+2}^{p,q}\cdot
j+g_{n+3}^{p,q}\cdot k,$ we have $\mathbf{t}\left( G_{n}^{p,q}\right)
=g_{n}^{p,q}$ and $\mathbf{n}\left( G_{n}^{p,q}\right) =g_{2n}^{a,b},$ with $%
a$ and $b$ as in Proposition 3.1. From here, we have that $\mathbf{n}\left(
G_{n}^{p,q}\text{-}\mathbf{t}\left( G_{n}^{p,q}\right) \right)
=g_{2n}^{a,b}-g_{n}^{p,q}=g_{2n}^{a,b}-\left( pf_{n-1}+ql_{n}\right) ^{2}.$
Using Proposition 2.1, relations v), ix) and xii), it results\newline
$\mathbf{n}\left( G_{n}^{p,q}\text{-}\mathbf{t}\left( G_{n}^{p,q}\right)
\right) $=$g_{2n}^{a,b}$-$\frac{p^{2}}{5}(l_{2n-2}$-$(-1)^{n}l_{0})$-$%
q^{2}\left( l_{2n}+2\left( -1\right) ^{n}\right) $-$2pq\left(
f_{2n-1}+\left( -1\right) ^{n}\right) =$\newline
$=g_{2n}^{a,b}-\frac{p^{2}}{5}l_{2n-2}-\frac{2p^{2}}{5}%
(-1)^{n}-q^{2}l_{2n}-2q^{2}\left( -1\right) ^{n}-2pqf_{2n-1}-2pq\left(
-1\right) ^{n}=$\newline
$=g_{2n}^{a,b}+g_{2n}^{-2pq,-q^{2}}-\frac{p^{2}}{5}l_{2n-2}-\frac{2p^{2}}{5}%
(-1)^{n}-2q^{2}\left( -1\right) ^{n}-2pq\left( -1\right) ^{n}=$\newline
$=g_{2n}^{a-2pq,b-q^{2}}-\frac{p^{2}}{5}l_{2n-2}-\frac{2p^{2}}{5}%
(-1)^{n}-2q^{2}\left( -1\right) ^{n}-2pq\left( -1\right) ^{n}=$\newline
$=g_{2n}^{a-2pq,b-q^{2}}-\frac{p^{2}}{5}l_{2n-2}-A,$ where $A=\frac{2p^{2}}{5%
}(-1)^{n}+2q^{2}\left( -1\right) ^{n}+2pq\left( -1\right) ^{n}.$ Using
relation $\left( 3.8\right) ,$ we obtain $x=2\lambda _{1}\left( G_{n}^{p,q}%
\text{-}g_{n}^{p,q}\right) +2\lambda _{2}\left( g_{2n}^{a-2pq,b-q^{2}}-\frac{%
p^{2}}{5}l_{2n-2}-A\right) =$\newline
$=G_{n}^{2\lambda _{1}p,2\lambda _{1}q}+g_{n}^{-2\lambda _{1}p,-2\lambda
_{1}q}+g_{2n}^{2\lambda _{2}a-4\lambda _{2}pq,2\lambda _{2}b-2\lambda
_{2}q^{2}}-2\lambda _{2}\frac{p^{2}}{5}l_{2n-2}-2\lambda _{2}A,\lambda
_{1},\lambda _{2}\in \mathbb{Q}.\square \medskip $

\textbf{Conclusions.} In this paper, starting from a special set of
elements, namely Fibonacci-Lucas quaternions, we proved that this set is an
order --in the sense of ring theory-- of the quaternion algebra $\mathbb{H}_{%
\mathbb{Q}}(\alpha ,\beta )$ and it can have an algebra structure over~$%
\mathbb{Q}.$ It will be very interesting to find all properties of this
algebra and to find conditions under which it is a division algebra or a
split algebra.

\textbf{Acknowledgments}

The authors thank Professor Rafa\l\ Ab\l amowicz and anonymous referees for
their comments, suggestions and ideas which helped us to improve this paper.

\begin{equation*}
\end{equation*}%
\textbf{References}\newline
\begin{equation*}
\end{equation*}%
[Ak, Ko, To; 14] M. Akyigit, HH Kosal, M. Tosun, \textit{\ Fibonacci
Generalized Quaternions }, Adv. Appl. Clifford Algebras, vol. \textbf{24},
issue: 3 (2014), p. 631-641\newline
[Fl, Sa; 14] C. Flaut, D. Savin, \textit{Some properties of the symbol
algebras of degree} $3$, Math. Reports, vol. \textbf{16(66)}(3)(2014),
p.443-463.\newline
[Fl, Sa, Io; 13] C. Flaut, D. Savin, G. Iorgulescu, \textit{Some properties
of Fibonacci and Lucas symbol elements}, Journal of Mathematical Sciences:
Advances and Applications, vol. \textbf{20} (2013), p. 37-43 \newline
[Fl, Sh; 13] C. Flaut, V. Shpakivskyi, \textit{On Generalized Fibonacci
Quaternions and Fibonacci-Narayana Quaternions,} accepted in Adv. Appl.
Clifford Algebras, vol. \textbf{23} issue 3 (2013), p. 673-688.\newline
[Fl, Sh; 13(1)] C. Flaut, V. Shpakivskyi, \textit{Real matrix
representations for the complex quaternions}, Adv. Appl. Clifford Algebras, 
\textbf{23(3)(2013)}, 657-671.\newline
[Fl; 06] C. Flaut, \textit{Divison algebras with dimension 2\symbol{94}t, t}$%
\in $\textit{N}, An. St. Univ. Ovidius Constanta, 13(2)(2006), 31-38.\newline
[Fl; 01] C. Flaut, \textit{Some equations in algebras obtained by the
Cayley-Dickson process}, An. St. Univ. Ovidius Constanta, \textbf{9(2)}%
(2001), 45-68.\newline
[Fl, Sa; 15] C. Flaut, D. Savin, \textit{Some examples of division symbol
algebras of degree 3 and 5}, accepted in Carpathian J Math. (arXiv:1310.1383
);\newline
[Fl, Sa; 14] C. Flaut, D. Savin, \textit{About quaternion algebras and
symbol algebras}, Bulletin of the Transilvania University of Brasov, \textbf{%
7(2)(56)(2014)}, 59-64.\newline
[Fl, St; 09] C. Flaut, M. \c{S}tef\~{a}nescu, \textit{Some equations over
generalized quaternion and octonion division algebras},  Bull. Math. Soc.
Sci. Math. Roumanie, 52(100), no. 4 (2009), 427-439.\newline
[Ha; 00] M. Hazewinkel, \textit{Handbook of Algebra}, Vol. 2, North Holland,
Amsterdam, 2000.\newline
[Ho; 61] A. F. Horadam, \textit{A Generalized Fibonacci Sequence}, Amer.
Math. Monthly, \textbf{68}(1961), 455-459.\newline
[Ho; 63] A. F. Horadam, \textit{Complex Fibonacci Numbers and Fibonacci
Quaternions}, Amer. Math. Monthly, \textbf{70}(1963), 289-291.\newline
[Lam; 04] T. Y. Lam, \textit{Introduction to Quadratic Forms over Fields,}
American Mathematical Society, 2004.\newline
[Sa; 14] D. Savin, \textit{About some split central simple algebras}, An.
Stiin. Univ. "`Ovidius" Constanta, Ser. Mat, \textbf{22} (1) (2014), p.
263-272.\newline
[Sa, Fl, Ci; 09] D. Savin, C. Flaut, C. Ciobanu, \textit{Some properties of
the symbol algebras}, Carpathian Journal of Mathematics, \textbf{25(2)(2009)}%
, p. 239-245\newline
[Vo; 10] J. Voight, The Arithmetic of Quaternion Algebras. Available on the
author's website: http://www.math.dartmouth.edu/ jvoight/ crmquat/book/quat-
modforms-041310.pdf, 2010.\newline
[Fib.] http://www.maths.surrey.ac.uk/hosted-sites/R.Knott/Fibonacci/fib.html%
\newline

\begin{equation*}
\end{equation*}

\bigskip

Cristina FLAUT

{\small Faculty of Mathematics and Computer Science, Ovidius University,}

{\small Bd. Mamaia 124, 900527, CONSTANTA, ROMANIA}

{\small http://cristinaflaut.wikispaces.com/;
http://www.univ-ovidius.ro/math/}

{\small e-mail: cflaut@univ-ovidius.ro; cristina\_flaut@yahoo.com}

\medskip \medskip \qquad\ \qquad\ \ 

Diana SAVIN

{\small Faculty of Mathematics and Computer Science, }

{\small Ovidius University, }

{\small Bd. Mamaia 124, 900527, CONSTANTA, ROMANIA }

{\small http://www.univ-ovidius.ro/math/}

{\small e-mail: \ savin.diana@univ-ovidius.ro, \ dianet72@yahoo.com}\bigskip
\bigskip

\end{document}